**Abstract**

An optimum solution free from degeneration is found to the system of linear algebraic equations with empirical coefficients and right-hand sides. The quadratic risk of estimators of the unknown solution vector is minimized over a class of linear systems with given square norm of the coefficient matrix and length of the right-hand side vector. Empirical coefficients and right-hand sides are assumed to be independent and normal with known variance. It is found that the optimal estimator has the form of a regularized minimum square solution with an extension multiple. A simple formula is derived showing explicitly the dependence of the minimal risk on parameters.

# 1 Introduction

The standard solution to systems of linear equations with random coefficients $R\mathbf{x} = \mathbf{y}$, where $R$ is the coefficient matrix, and $\mathbf{x}$ is a vector of unknowns can be unstable or may not exist if the variance of coefficients is sufficiently large. This can be caused by an incorrect solution and by the inconsistency of the system of random equations. The regularization by A.N.Tikhonov [1] implies an artificial requirement of minimal complexity; this procedure guaranties some solution, but this solution provides neither the minimum of the quadratic risk, nor the minimum of residuals.

The estimation of unknowns over a single observation of $R$ and $\mathbf{y}$ by the minimum square method leads to the solution $\mathbf{x} = (\mathbf{R}^{\mathbf{T}}\mathbf{R})^{-1}\mathbf{R}^{\mathbf{T}}\mathbf{y}$, which also can be unstable or non-existing. The confluent analysis methods [2] provide solutions of the form $\mathbf{x} = (\mathbf{R}^{\mathbf{T}}\mathbf{R} - \lambda \mathbf{I})^{-1}\mathbf{R}^{\mathbf{T}}\mathbf{y}$, where $\lambda \geq 0$, and $I$ is the identity matrix. These estimates are even more sensible to the coefficients inaccuracy and certainly do not exist when the minimum square solution does not exist (that is related to the necessity of a simultaneous additional estimation of the matrix coefficients). The extremum problem that is solved in the confluent analysis is based on the application of the maximum likelihood estimations. However, it is well known [3] that the maximum likelihood estimators can have the quadratic risk substantially greater than the minimal one if the number of parameters is comparable with the sample size.

One can expect that the a priori information on the distribution of random coefficients and right-hand sides can be used to improve solutions, and even by a single observation of the coefficient matrix if the systems of equations is large and random inaccuracies of the coefficients are independent.

In paper [3] and monograph [4] by V.L.Girko, the problem of solutions to large systems of random equations is investigated. It is proved that the standard minimum square method using the inversion of empirical coefficient matrix leads to solutions with a substantial bias. The author suggests asymptotically unbiased estimators of the unknowns vector (G-estimators). However, the quadratic risk of these unbiased estimators also proves to be not minimal.

In [5], the asymptotic approach of [3] was applied to find an asymptotically unimprovable linear operator solving systems of random linear equations of increa-sing dimension. However, the possibilities to decrease the quadratic risk of solutions for finite systems were not yet investigated.

In this paper, we obtain an estimator of the solution vector minimizing the square risk on the average for a class of systems with arbitrary fixed numbers of equations and unknowns.

## 2 Optimal solution

Let us find a general form of an unimprovable solution in the Bayes approach by averaging over possible equations and their solutions.

Suppose a consistent system of linear equations is given $A\mathbf{x} = \mathbf{b}$, where $A$ is the coefficient matrix of $N$ rows and $n$ columns, $N \geq n$, $\mathbf{x}$ is the unknown vector of the dimension $n$, and $\mathbf{b}$ is the right-hand side vector of the dimension $N$. The observer only knows the matrix $R = A + \delta A$, and vector $\mathbf{y} = \mathbf{b} + \delta \mathbf{b}$, where $\delta A$ are matrices with random entries $\delta A_{ij} \sim \mathbf{N}(0, p/n)$ and $\delta \mathbf{b}$ are vectors with components $\delta b_i \sim \mathbf{N}(0, q/n)$. Assume that all these random values are independent, $i = 1, \ldots, N, \quad j = 1, \ldots, n$.

We construct a solution procedure that is unimprovable on the average for a set of problems with different matrices $A$ and vectors $\mathbf{b}$, and consider Bayes distributions $A_{ij} \sim \mathbf{N}(0, a/n)$ and $\mathbf{x_j} \sim \mathbf{N}(0, 1/n)$ for all $i$ and $j$ assuming all random values independent. It corresponds to a uniform distribution of entries of $A$ given the quadratic norm distributed as $\chi^2$ and a uniform distribution of directions of vectors $\mathbf{x}$. We minimize the quadratic risk of the estimator $\widehat{\mathbf{x}}$

$$D = D(\widehat{\mathbf{x}}) = \mathbf{E}\ (\mathbf{x} - \widehat{\mathbf{x}})^2. \tag{1}$$

Here (and in the following) the square of a vector denotes the square of its length, and the expectation is calculated over the joint distribution of $A, \delta A, \mathbf{x}$, and $\delta \mathbf{b}$.

Let us restrict ourselves with a class K of regularized pseudo solutions of the form

$$\widehat{\mathbf{x}} = \Gamma R^T \mathbf{y}, \tag{2}$$

where the matrix $\Gamma = \Gamma(R^T R)$ is diagonalized together with $R^T R$ and has eigenva-lues $\gamma(\lambda)$ corresponding to the eigenvalues $\lambda$ of $R^T R$. Restrict functions $\gamma(u)$ by the requirement that it is a non-negative measurable function such that the product $u(1+u)\gamma^2(u)$ is bounded by a constant.

Note that, for the standard solution, $\gamma(u) = 1/u \in \mathrm{K}$.

Define a non-random distribution function for the set of eigenvalues of the empirical matrices $R^T R$ as the expectation

$$F(u) = \mathbf{E}\ \mathbf{n}^{-1} \sum_{\mathbf{i=1}}^{\mathbf{n}} \mathrm{ind}(\lambda_\mathbf{i} \leq \mathbf{u}),$$

where $\lambda_i$ are eigenvalues of $R^T R$. Using this function we can write, for example,

$$\mathbf{E}\ n^{-1}\mathrm{tr}(\mathbf{R}\Gamma^2\mathbf{R}^T) = n^{-1}\mathbf{E}\sum_{i=1}^{n}\lambda_i \gamma^2(\lambda_i) = \int u\gamma^2(u)dF(u).$$

THEOREM 1 Let an estimator $\widehat{\mathbf{x}} \in K$ of the solution 2 to the system $A\mathbf{x} = \mathbf{b}$ be calculated using the observed empiric matrix $R = A + \delta A$ and vector $\mathbf{y} = \mathbf{b} + \delta\mathbf{b}$, where entries of the matrices $A = \{A_{ij}\}$ and $\delta A = \{\delta A_{ij}\}$ and components of the vectors $\mathbf{x} = \{x_j\}$ and $\delta\mathbf{b} = \{\delta b_i\}$ are independent random values

$$A_{ij} \sim N(0, a/n), \quad x_j \sim N(0, 1/n)), \quad \delta A_{ij} \sim N(0, p/n), \quad \delta b_i \sim N(0, q/n),$$

$i = 1, \ldots, N, \quad j = 1, \ldots, n$, where $a > 0$. Then the functional 1 equals

$$D(\widehat{\mathbf{x}}) = \int \left[1 - 2\theta u\gamma(u) + \theta^2 u^2 \gamma^2(u) + su\gamma^2(u)\right]\,dF(u), \tag{3}$$

where

$$\theta = \frac{a}{a+p}, \quad s = \frac{ap}{a+p} + q.$$

Note that the expression 3 is quadratic with respect to $\gamma(u)$ and allows the standard minimization.

COROLLARY Let $a > 0$. The expression 3 reaches the minimum in the class K for $\gamma(u) = \gamma_{opt}(u) = \theta/(\theta^2 u + s)$,

$$\widehat{\mathbf{x}} = \widehat{\mathbf{x}}_{opt} = \frac{\theta}{\theta^2 R^T R + sI}\ R^T \mathbf{y},$$

and the minimal value of 1 is

$$D = D(\widehat{\mathbf{x}}_{opt}) = \int \frac{s}{s + \theta^2 u}dF(u).$$

## 3 Discussion

Thus we can draw a conclusion that the solution of empirical linear algebraic equations with normal errors can be stabilized and made more accurate by replacing standard minimum square solution by an optimal on the average regular estimator

$$\widehat{\mathbf{x}} = \alpha(R^T R + tI)^{-1} R^T \mathbf{y} \tag{4}$$

with a "ridge"-parameter $t = s/\theta^2$ and a scalar scaling coefficient $\alpha = 1/\theta$. Here $\alpha \geq 1$ presents an *extension* coefficient in contrast to shrinkage multiples for well-known Stein estimators. The matrix $\alpha(R^T R + tI)^{-1}$ in 4 is, actually, a scaled and regularized estimator of $(A^T A)^{-1}$ optimal with respect to the minimization of $D$. Such solutions are characteristic for a number of extremum problems solved in the theory of high-dimensional multivariate statistical analysis [6].

We can compare the optimal value $D_{\text{opt}}$ with the square risk of the standard solution $D_{\text{std}}$ if we set $\gamma(u) = 1/u$ in 3,

$$D_{std} = \int [(1-\theta)^2 + \frac{s}{u}] \, dF(u). \tag{5}$$

It is easy to prove that the difference $D_{std} - D_{\text{opt}} \geq 0$ and equals 0 for $q = 0$ and $p = 0$.

For $p = 0$ and $q > 0$, we have $\theta = 1$, the matrix $R$ is non-random, $\alpha = 1$, $t = q$ and $D_{\text{opt}}$ is less $D_{\text{std}}$ due to the factor $u/(u+q)$ in the integrals over $dF(u)$.

If $p > 0$ or $q > 0$, the expression 5 contains $\int u^{-1} dF(u)$ with a non-zero coefficient. Note that matrices $W = R^T R$ are the Wishart $n \times n$ matrices, calculated for variables distributed as $\mathbf{N}(\mathbf{0}, (\mathbf{a}+\mathbf{p})\mathbf{N/n})$ over a sample of size $N$. The integral $\int u^{-1} dF(u) = \mathbf{E} \, \mathbf{n}^{-1} \text{tr} \mathbf{W}^{-1}$, where the right-hand is infinite for $n \geq N$. Thus, for $n = N$ the quadratic risk of standard solution is infinitely large while $D_{\text{opt}} = 1$ for the optimal solution.

## Proof of the theorem 1.

Note that $D$ is a sum of three addends

$$D = \mathbf{E} \, \mathbf{x^T x} - 2\mathbf{E} \, \mathbf{x^T \Gamma R^T y} + \mathbf{E} \, \mathbf{y^T R \Gamma^2 R^T y}.$$

Let angular brackets denote the normed trace of a matrix: $\langle M \rangle = n^{-1} \text{tr} M$.

First, we average over the distribution of $\delta \mathbf{b}$. Let us use the obvious property of normal distributions

$$\mathbf{E} \, (\delta \mathbf{b^T} \mathbf{M} \delta \mathbf{b}) = \mathbf{q} \mathbf{n^{-1}} \text{tr} \mathbf{M},$$

where $M$ is an $N \times N$ matrix of constants. We find that

$$D = \mathbf{E} \, \mathbf{x^T x} - 2\mathbf{E} \, \mathbf{x^T \Gamma R^T b} + \mathbf{E} \, \mathbf{b^T R \Gamma^2 R^T b} + \mathbf{qE} \, \langle \mathbf{R \Gamma^2 R^T} \rangle.$$

We average with respect to the distribution of $\mathbf{x}$. In view of the equality $\mathbf{b} = \mathbf{A}\mathbf{x}$ we have

$$D = 1 - 2\mathbf{E} \, \langle \mathbf{\Gamma R^T A} \rangle + \mathbf{E} \, \langle \mathbf{A^T R \Gamma^2 R^T A} \rangle + \mathbf{qE} \, \langle \mathbf{R\Gamma^2 R^T} \rangle. \tag{6}$$

In the further transformations, we use a simple property of normal distributions: if a random value $r$ is normally distributed with zero average and the variance

$\sigma^2$, then for any differentiable function $f(\cdot)$, we have

$$\mathbf{E}\ \mathbf{r} f(\mathbf{r}) = \sigma^2 \mathbf{E}\ \frac{\partial \mathbf{f}}{\partial \mathbf{r}}. \tag{7}$$

For example, consider the second term in 3. Apply 7 to the normal variable $A$. We obtain

$$\mathbf{E}\ \langle \mathbf{\Gamma R^T A} \rangle = \mathbf{n}^{-1} \mathbf{E}\ \mathbf{A_{ij}}(\mathbf{\Gamma R^T})_{\mathbf{ji}} = \frac{\mathbf{a}}{\mathbf{n^2}} \mathbf{E}\ \frac{\partial(\mathbf{\Gamma R^T})_{\mathbf{ji}}}{\partial \mathbf{A_{ij}}},$$

where (and in the following) the summation over repeated indexes is implied.

Note that all entries of the random matrix $R$ have the joint variance $(a+p)/n$. Therefore,

$$\mathbf{E}\ \langle \mathbf{\Gamma R^T R} \rangle = \mathbf{n}^{-1} \mathbf{E}\ \mathbf{R_{ij}}(\mathbf{\Gamma R^T})_{\mathbf{ji}} = \frac{\mathbf{a+p}}{\mathbf{n^2}}\ \frac{\partial(\mathbf{\Gamma R^T})_{\mathbf{ji}}}{\partial \mathbf{R_{ij}}}.$$

Let $a > 0$. Comparing these two expressions we find that

$$\mathbf{E}\ \langle \mathbf{\Gamma R^T A} \rangle = \mathbf{a(a+p)}^{-1} \mathbf{E}\ \langle \mathbf{\Gamma R^T R} \rangle$$

Similarly, one can derive the relation

$$\mathbf{E}\ \langle \mathbf{R^T R \Gamma^2 R^T A} \rangle = \mathbf{a(a+p)}^{-1} \mathbf{E}\ \langle \mathbf{R^T R \Gamma^2 R^T R} \rangle.$$

Once again, we use 7 with respect to the random $A$ and obtain:

$$\mathbf{E}\ \langle \mathbf{A^T R \Gamma^2 R^T A} \rangle = \frac{\mathbf{a}}{\mathbf{n^2}} \mathbf{E}\ \mathbf{A_{ij}} \frac{\partial(\mathbf{R \Gamma^2 R^T})_{\mathbf{jk}}}{\partial \mathbf{A_{kj}}} + \mathbf{aE}\ \langle \mathbf{R \Gamma^2 R^T} \rangle,$$

where (and in the following) $k = 1, \ldots, N$. Further, we use 7 with respect to $\delta A$ and obtain

$$\mathbf{E}\ \langle \mathbf{A^T R \Gamma^2 R^T} \delta \mathbf{A} \rangle = \frac{\mathbf{p}}{\mathbf{n^2}} \mathbf{E}\ \mathbf{A_{ji}} \frac{\partial(\mathbf{R \Gamma^2 R^T})_{\mathbf{jk}}}{\partial (\delta \mathbf{A})_{\mathbf{ki}}}.$$

Adding two latter equalities we have

$$\mathbf{E}\ \langle \mathbf{A^T R \Gamma^2 R^T R} \rangle = \frac{\mathbf{a+p}}{\mathbf{n^2}} \mathbf{E}\ \mathbf{A_{ji}} \frac{\partial(\mathbf{R \Gamma^2 R^T})_{\mathbf{jk}}}{\partial \mathbf{A_{ki}}} + \mathbf{aE}\ \langle \mathbf{R \Gamma^2 R^T} \rangle.$$

Now we combine the five latter equalities and obtain

$$\mathbf{E}\ \langle \mathbf{A^T R \Gamma^2 R^T A} \rangle = \frac{\mathbf{a^2}}{(\mathbf{a+p})^2} \mathbf{E}\ \langle \mathbf{R^T R \Gamma^2 R^T R} \rangle + \frac{\mathbf{ap}}{\mathbf{a+p}} \mathbf{E}\ \langle \mathbf{R \Gamma^2 R^T} \rangle. \tag{8}$$

These relations are sufficient to express $D$ in terms of the observed variables only. Substitute 8 to 6. We have

$$D = 1 - 2\frac{a}{a+p}\mathbf{E}\,\langle \mathbf{\Gamma R^T R}\rangle + \frac{\mathbf{a^2}}{\mathbf{(a+p)^2}}\mathbf{E}\,\langle \mathbf{R^T R \Gamma^2 R^T R}\rangle + \mathbf{sE}\,\langle \mathbf{R\Gamma^2 R^T}\rangle. \quad (9)$$

Passing to the (random) system of coordinates in which $R^T R$ is diagonal we can replace formally

$$n^{-1}\sum_{i=1}^{n}\varphi(\lambda_i) = \int \varphi(u)dF(u),$$

where $\lambda_i$ are eigenvalues of $R^T R$. The expression 3 follows from 9. The proof of Theorem 1 is complete.